\documentclass[fleqn, eqnumis, pdf]{hjms}

\usepackage{footnote}
\usepackage{multicol}
\usepackage{booktabs}
\usepackage[flushleft]{threeparttable}
\usepackage[font=small,labelfont=bf]{caption}

\begin{document}

\title{Existence of a unique positive solution for a singular fractional boundary value problem}

\author{ Karimov E.T.\footnote{ Department of Mathematics and Statistics, Sultan Qaboos University, Al-Khoud 123, Muscat, Oman,  \ Email: {\tt erkinjon@gmail.com }}\footnote{Corresponding author} and 
Sadarangani K.\footnote{ Departamento de Matematicas, Universidad de Las Palmas de Gran Canaria, Campus de Tafira Baja, 35017 Las Palmas de Gran Canaria, Spain,  \ Email: {\tt ksadaran@dma.ulpgc.es }}}
{ }

\begin{abstract}
In the present work, we discuss the existence of a unique positive solution of a boundary value problem for nonlinear fractional order equation with singularity. Precisely, order of equation $D_{0+}^\alpha u(t)=f(t,u(t))$ belongs to $(3,4]$ and $f$ has a singularity at $t=0$ and as a boundary conditions we use $u(0)=u(1)=u'(0)=u'(1)=0$. Using fixed point theorem, we prove the existence of unique positive solution of the considered problem.

\end{abstract}



\section{Introduction}

In this paper, we study the existence and uniqueness of positive solution for the following singular fractional boundary value problem
$$
\left\{
\begin{array}{l}
D_{0+}^\alpha u(t)=f(t,u(t)),\,\,0<t<1\\
u(0)=u(1)=u'(0)=u'(1)=0,
\end{array}
\right.
\eqno (1)
$$
where $\alpha\in (3,4]$, and $D_{0+}^\alpha$ denotes the Riemann-Liouville fractional derivative. Moreover, $f: (0,1]\times [0,\infty)\rightarrow [0,\infty)$ with $\lim\limits_{t\rightarrow 0+} f(t,-)=\infty$ (i.e. $f$ is singular at $t=0$).

Similar problem was investigated in [1], in case when $\alpha\in (1,2]$ and with boundary conditions $u(0)=u(1)=0$. We note as well work [2], where
the following problem
$$
\left\{
\begin{array}{l}
D^\alpha u+f(t,u,u',D^\mu u)=0,\,\,0<t<1\\
u(0)==u'(0)=u'(1)=0,
\end{array}
\right.
$$
was under consideration. Here $\alpha\in (2,3), \mu\in (0,1)$ and function $f(t,x,y,z)$ is singular at the value of $0$ of its arguments $x,y,z$.

We would like notice some related recent works [3-5], which consider higher order fractional nonlinear equations for the subject of the existence of positive solutions.

\section{Preliminaries}

We need the following lemma, which appear in [6].

{\bf Lemma 1.} (Lemma 2.3 of [6]) Given $h\in C[0,1]$ and $3<\alpha\leq 4$, a unique solution of
$$
\left\{
\begin{array}{l}
D_{0+}^\alpha u(t)=h(t),\,\,0<t<1\\
u(0)=u(1)=u'(0)=u'(1)=0,
\end{array}
\right.\eqno (1)
$$
is
$$
u(t)=\int\limits_0^1G(t,s)h(s)ds,
$$
where
$$
G(t,s)=
\left\{
\begin{array}{l}
\frac{(t-s)^{\alpha-1}+(1-s)^{\alpha-2}t^{\alpha-2}\left[(s-t)+(\alpha-2)(1-t)s\right]}{\Gamma(\alpha)},\,\,0\leq s\leq 1\\
\frac{(1-s)^{\alpha-2}t^{\alpha-2}\left[(s-t)+(\alpha-2)(1-t)s\right]}{\Gamma(\alpha)},\,\,0\leq t\leq s\leq 1\\
\end{array}
\right.
$$

{\bf Lemma 2.} (Lemma 2.4 of [6]) The function $G(t,s)$ appearing in Lemma 1 satisfies:

(a) $G(t,s)>0$ for $t,s \in (0,1)$;

(b) $G(t,s)$ is continuous on $[0,1]\times [0,1]$.

For our study, we need a fixed point theorem. This theorem uses the following class of functions $\mathfrak{F}$.

By $\mathfrak{F}$ we denote the class of functions $\varphi: (0, \infty)\rightarrow \mathbb{R}$ satisfying the following conditions:

(a) $\varphi$ is strictly increasing;

(b) For each sequence $(t_n)\subset (0,\infty)$
$$
\lim\limits_{n\rightarrow \infty} t_n=0  \,\,\Leftrightarrow \lim\limits_{n\rightarrow \infty}\varphi(t_n)=-\infty;
$$

(c) There exists $k\in (0,1)$ such that $\lim\limits_{t\rightarrow 0^+}t^k\varphi(t)=0$.

Examples of functions belonging to $\mathfrak{F}$ are $\varphi(t)=-\frac{1}{\sqrt{t}}, \varphi(t)=\ln t, \varphi(t)=\ln t+t,\varphi(t)=\ln (t^2+t)$.

The result about fixed point which we use is the following and it appears in [7]:

\textbf{Theorem 3.} Let $(X,d)$ be a complete metric space and $T:X\rightarrow X$ a mapping such that there exist $\tau>0$ and $\varphi\in \mathfrak{F}$ satisfying for any $x,y\in X$ with $d(Tx,Ty)>0$,
$$
\tau+\varphi\left(d(Tx,Ty)\right)\leq \varphi(d(x,y)).
$$
Then $T$ has a unique fixed point.

\section{Main result}

Our starting point of this section is the following lemma.

\textbf{Lemma 4.} Let $0<\sigma<1$, $3<\alpha<4$ and $F:(0,1]\rightarrow \mathbb{R}$ is continuous function with $\lim\limits_{t\rightarrow 0^+} F(t)=\infty$. Suppose that $t^\sigma F(t)$ is a continuous function on $[0,1]$. Then the function defined by
$$
H(t)=\int\limits_0^1 G(t,s)F(s)ds
$$
is continuous on $[0,1]$, where $G(t,s)$ is the Green function appearing in Lemma 1.

\emph{Proof:} We consider three cases:

{\bf \texttt {1. Case No 1}}. $t_0=0$.

It is clear that $H(0)=0$. Since $t^\sigma F(t)$ is continuous on $[0,1]$, we can find a constant $M>0$ such that
$$
|t^\sigma F(t)|\leq M \,\,\,for \,\,any\,\,\, t\in [0,1].
$$
Moreover, we have
$$
\begin{array}{l}
|H(t)-H(0)|=|H(t)|=\left| \int\limits_0^1 G(t,s)F(s)ds\right|=\left| \int\limits_0^1 G(t,s)s^{-\sigma} s^\sigma F(s)ds\right|=\\
=\left| \int\limits_0^t \frac{(t-s)^{\alpha-1}+(1-s)^{\alpha-2}t^{\alpha-2}[(s-t)+(\alpha-2)(1-t)s]}{\Gamma(\alpha)}s^{-\sigma}s^\sigma F(s)ds\right.\\
\left.+
\int\limits_t^1 \frac{(1-s)^{\alpha-2}t^{\alpha-2}[(s-t)+(\alpha-2)(1-t)s]}{\Gamma(\alpha)}s^{-\sigma}s^\sigma F(s)ds\right|=\\
=\left|\int\limits_0^1 \frac{(1-s)^{\alpha-2}t^{\alpha-2}[(s-t)+(\alpha-2)(1-t)s]}{\Gamma(\alpha)}s^{-\sigma}s^\sigma F(s)ds+
\int\limits_0^t \frac{(t-s)^{\alpha-1}}{\Gamma(\alpha)}s^{-\sigma}s^\sigma F(s)ds\right|\leq\\
\leq \frac{Mt^{\alpha-2}}{\Gamma(\alpha)}\int\limits_0^1 (1-s)^{\alpha-2}|(s-t)+(\alpha-2)(1-t)s|s^{-\sigma}ds+
\frac{M}{\Gamma(\alpha)}\int\limits_0^t (t-s)^{\alpha-1}s^{-\sigma}ds\leq\\
\leq \frac{M(\alpha-1)t^{\alpha-2}}{\Gamma(\alpha)}\int\limits_0^1 (1-s)^{\alpha-2}s^{-\sigma}ds+
\frac{Mt^{\alpha-1}}{\Gamma(\alpha)}\int\limits_0^t \left(1-\frac{s}{t}\right)^{\alpha-1}s^{-\sigma}ds.\\
\end{array}
$$
Considering definition of Euler's beta-function, we derive
$$
|H(t)-H(0)|\leq \frac{M(\alpha-1)t^{\alpha-2}}{\Gamma(\alpha)} B(1-\sigma,\alpha-1)+\frac{Mt^{\alpha-\sigma}}{\Gamma(\alpha)} B(1-\sigma, \alpha).
$$
From this we deduce that $|H(t)-H(0)|\rightarrow 0$ when $t\rightarrow 0$.

This proves that $H$ is continuous at $t_0=0$.

{\bf \texttt {2. Case No 2}}. $t_0\in (0,1)$.

We take $t_n\rightarrow t_0$ and we have to prove that $H(t_n)\rightarrow H(t_0)$. Without loss of generality, we consider $t_n>t_0$. Then, we have
$$
\begin{array}{l}
\left|H(t_n)-H(t_0)\right|=\left|\int\limits_0^{t_n}\frac{(t_n-s)^{\alpha-1}+(1-s)^{\alpha-2}t_n^{\alpha-2}[(s-t_n)+(\alpha-2)(1-t_n)s]}{\Gamma(\alpha)}
s^{-\sigma} s^\sigma F(s)ds+\right.\\
\left.+\int\limits_{t_n}^1\frac{(1-s)^{\alpha-2}t_n^{\alpha-2}[(s-t_n)+(\alpha-2)(1-t_n)s]}{\Gamma(\alpha)}
s^{-\sigma} s^\sigma F(s)ds-\right.\\
-\int\limits_0^{t_0}\frac{(t_0-s)^{\alpha-1}+(1-s)^{\alpha-2}t_0^{\alpha-2}[(s-t_0)+(\alpha-2)(1-t_0)s]}{\Gamma(\alpha)}
s^{-\sigma} s^\sigma F(s)ds-\\
\left.-\int\limits_{t_0}^1\frac{(1-s)^{\alpha-2}t_0^{\alpha-2}[(s-t_0)+(\alpha-2)(1-t_0)s]}{\Gamma(\alpha)}
s^{-\sigma} s^\sigma F(s)ds\right|=\\
=\left|\int\limits_0^1\frac{(1-s)^{\alpha-2}t_n^{\alpha-2}[(s-t_n)+(\alpha-2)(1-t_n)s]}{\Gamma(\alpha)}
s^{-\sigma} s^\sigma F(s)ds+\int\limits_0^{t_n}\frac{(t_n-s)^{\alpha-1}}{\Gamma(\alpha)}s^{-\sigma} s^\sigma F(s)ds-\right.\\
\left.-\int\limits_0^1\frac{(1-s)^{\alpha-2}t_0^{\alpha-2}[(s-t_0)+(\alpha-2)(1-t_0)s]}{\Gamma(\alpha)}
s^{-\sigma} s^\sigma F(s)ds-\int\limits_0^{t_0}\frac{(t_0-s)^{\alpha-1}}{\Gamma(\alpha)}s^{-\sigma} s^\sigma F(s)ds\right|=\\

=\left|\int\limits_0^1\frac{(1-s)^{\alpha-2}\left(t_n^{\alpha-2}-t_0^{\alpha-2}\right)[(s-t_n)+(\alpha-2)(1-t_n)s]}{\Gamma(\alpha)}
s^{-\sigma} s^\sigma F(s)ds+\right.\\
\left.+\int\limits_0^1\frac{(1-s)^{\alpha-2}t_0^{\alpha-2}\left[(s-t_n)+(\alpha-2)(1-t_n)s-[(s-t_0)+(\alpha-2)(1-t_0)s]\right]}{\Gamma(\alpha)}
s^{-\sigma} s^\sigma F(s)ds+\right.\\
\left.+\int\limits_0^{t_0}\frac{\left[(t_n-s)^{\alpha-1}-(t_0-s)^{\alpha-1}\right]}{\Gamma(\alpha)}s^{-\sigma} s^\sigma F(s)ds+\int\limits_{t_0}^{t_n}\frac{(t_n-s)^{\alpha-1}}{\Gamma(\alpha)}s^{-\sigma} s^\sigma F(s)ds\right|\leq\\
\leq \frac{M\left|t_n^{\alpha-2}-t_0^{\alpha-2}\right|}{\Gamma(\alpha)}(\alpha-1)\int\limits_0^1(1-s)^{\alpha-2}s^{-\sigma}ds+\\
+\frac{M t_0^{\alpha-2}}{\Gamma(\alpha)}\int\limits_0^1(1-s)^{\alpha-2}|t_n-t_0|(\alpha-1)s^{-\sigma}ds+\\
+\frac{M}{\Gamma(\alpha)}\int\limits_0^{t_0}\left|(t_n-s)^{\alpha-1}-(t_0-s)^{\alpha-1}\right|s^{-\sigma}ds+
\frac{M}{\Gamma(\alpha)}\int\limits_{t_0}^{t_n}(t_n-s)^{\alpha-1}s^{-\sigma}ds\leq\\
\leq \frac{M}{\Gamma(\alpha)}\left(t_n^{\alpha-2}-t_0^{\alpha-2}\right)(\alpha-1)B(1-\sigma,\alpha-1)+
\frac{M(t_n-t_0)}{\Gamma(\alpha)}(\alpha-1)B(1-\sigma,\alpha-1)+\\
+\frac{M}{\Gamma(\alpha)}I_n^1+\frac{M}{\Gamma(\alpha)}I_n^2,\\
\end{array}
$$
where
$$
I_n^1=\int\limits_0^{t_0}\left[(t_n-s)^{\alpha-1}-(t_0-s)^{\alpha-1}\right]s^{-\sigma}ds,\,\,\,I_n^2=\int\limits_{t_0}^{t_n}(t_n-s)^{\alpha-1}s^{-\sigma}ds.
$$

In the sequel, we will prove that $I_n^1\rightarrow 0$ when $n\rightarrow \infty$. In fact, as
$$
\left[(t_n-s)^{\alpha-1}-(t_0-s)^{\alpha-1}\right]s^{-\sigma}\leq \left[|t_n-s|^{\alpha-1}-|t_0-s|^{\alpha-1}\right]s^{-\sigma}\leq 2s^{-\sigma}
$$
and $\int\limits_0^1 2s^{-\sigma}ds=\frac{2}{1-\sigma}<\infty$. By Lebesque's dominated convergence theorem $I_n^1\rightarrow 0$ when $n\rightarrow \infty$.

Now, we will prove that $I_n^2\rightarrow 0$ when $n\rightarrow \infty$. In fact, since
$$
I_n^2=\int\limits_{t_0}^{t_n}(t_n-s)^{\alpha-1}s^{-\sigma}ds\leq \int\limits_{t_0}^{t_n}s^{-\sigma}ds=\frac{1}{1-\sigma}\left(t_n^{1-\sigma}-t_0^{1-\sigma}\right)
$$
and as $t_n\rightarrow t_0$, we obtain the desired result.

Finally, taking into account above obtained estimates, we infer that $|H(t_n)-H(t_0)|\rightarrow 0$ when $n\rightarrow \infty$.

{\bf \texttt {3. Case No 3}}. $t_0=1$.

It is clear that $H(1)=0$ and following the same argument that in Case No 1, we can prove that continuity of $H$ at $t_0=1$.

\textbf{Lemma 5.} Suppose that $0<\sigma<1$. Then there exists
$$
N=\max\limits_{0\leq t\leq 1}\int\limits_0^1 G(t,s)s^{-\sigma}ds.
$$

\textbf{\emph{Proof:}} Considering representation of the function $G(t,s)$ and evaluations of Lemma 4, we derive
$$
\begin{array}{l}
\int\limits_0^1 G(t,s)s^{-\sigma}ds=\frac{1}{\Gamma(\alpha)}\left[t^{\alpha-\sigma}B(1-\sigma, \alpha)-t^{\alpha-1}
\left(B(1-\sigma, \alpha-1)+\right.\right.\\
\left.\left.+(\alpha-2)B(2-\sigma, \alpha-1)\right)+(\alpha-1)t^{\alpha-2}B(2-\sigma, \alpha-1)\right].
\end{array}
$$
Taking
$$
B(1-\sigma,\alpha)=\frac{\alpha-1}{\alpha-\sigma}B(1-\sigma, \alpha-1);\,\,B(2-\sigma, \alpha-1)=\frac{1-\sigma}{\alpha-\sigma}B(1-\sigma, \alpha-1),
$$
into account we infer
$$
\begin{array}{l}

\int\limits_0^1 G(t,s)s^{-\sigma}ds=\\
=\frac{B(1-\sigma, \alpha-1)}{\Gamma(\alpha)}\left[\frac{\alpha-1}{\alpha-\sigma}t^{\alpha-\sigma}-\left(1+\frac{(\alpha-2)(1-\sigma)}{\alpha-\sigma}\right)t^{\alpha-1}
+\frac{(\alpha-1)(1-\sigma)}{\alpha-\sigma}t^{\alpha-2}\right].
\end{array}
$$
Denoting $L(t)=\int\limits_0^1 G(t,s)s^{-\sigma}ds$, from the last equality one can easily derive that $L(0)=0,\, L(1)=0$. Since $G(t,s)\geq 0$, then $L(t)\geq 0$ and as $L(t)$ is continuous on $[0,1]$, it has a maximum. This proves Lemma 5.

\textbf{Theorem 6.} Let $0<\sigma<1, \, 3<\alpha\leq 4,\, f: \, (0,1]\times [0,\infty)$ be continuous and $\lim\limits_{t\rightarrow 0^+} f(t,\cdot)=\infty$, $t^\sigma f(t,y)$ be continuous function on $[0,1]\times [0,\infty)$. Assume that there exist constants $0<\lambda\leq \frac{1}{N}$, and $\tau>0$ such that for $x,y \in [0,\infty)$ and $t\in [0,1]$
$$
t^\sigma |f(t,x)-f(t,y)|\leq \frac{\lambda |x-y|}{\left(1+\tau\sqrt{|x-y|}\right)^2}.
$$
Then Problem (1) has a unique non-negative solution.

\textbf{\emph{Proof:}} Consider the cone $P=\left\{u\in C[0,1]:\, u\geq 0\right\}$. Notice that $P$ is a closed subset of $C[0,1]$ and therefore, $(P,d)$ is a complete metric space where
$$
d(x,y)=\sup\left\{|x(t)-y(t)|:\, t\in[0,1]\right\}\,\, for\,\, x,y\in P.
$$

Now, for $u\in P$ we define the operator $T$ by
$$
(Tu)(t)=\int\limits_0^1 G(t,s) f(s,u(s))ds=\int\limits_0^1 G(t,s) s^{-\sigma} s^\sigma f(s,u(s))ds.
$$

In virtue of Lemma 4, for $u\in P$, $Tu\in C[0,1]$ and, since $G(t,s)$ and $t^\sigma f(t,y)$ are non-negative functions, $Tu\geq 0$ for $u\in P$. Therefore, $T$ applies $P$ into itself.

Next, we check that assumptions of Theorem 3 are satisfied. In fact, for $u,v\in P$ with $d(Tu,Tv)>0$, we have
$$
\begin{array}{l}
d(Tu,Tv)=\max\limits_{t\in [0,1]} \left|(Tu)(t)-(Tv)(t)\right|=\\
=\max\limits_{t\in [0,1]}\left|\int\limits_0^1 G(t,s)s^{-\sigma}s^\sigma \left(f(s,u(s))-f(s,v(s))\right)ds\right|\leq\\
\leq \max\limits_{t\in [0,1]}\int\limits_0^1 G(t,s)s^{-\sigma}s^\sigma \left|f(s,u(s))-f(s,v(s))\right|ds\leq\\
\leq \max\limits_{t\in [0,1]}\int\limits_0^1 G(t,s)s^{-\sigma}\frac{\lambda|u(s)-v(s)|}{\left(1+\tau\sqrt{|u(s)-v(s)|}\right)^2}ds\leq
\max\limits_{t\in [0,1]}\int\limits_0^1 G(t,s)s^{-\sigma}\frac{\lambda d(u,v)}{\left(1+\tau\sqrt{d(u,v)}\right)^2}ds=\\
=\frac{\lambda d(u,v)}{\left(1+\tau\sqrt{d(u,v)}\right)^2}\max\limits_{t\in [0,1]}\int\limits_0^1 G(t,s)s^{-\sigma}ds=
\frac{\lambda d(u,v)}{\left(1+\tau\sqrt{d(u,v)}\right)^2} N\leq \frac{d(u,v)}{\left(1+\tau\sqrt{d(u,v)}\right)^2},\\
\end{array}
$$
where we have used that $\lambda\leq\frac{1}{N}$ and the non-decreasing character of the function $\beta(t)=\frac{t}{\left(1+\tau\sqrt{t}\right)^2}$.
Therefore,
$$
d(Tu,Tv)\leq \frac{d(u,v)}{\left(1+\tau\sqrt{d(u,v)}\right)^2}.
$$
This gives us
$$
\sqrt{d(Tu,Tv)}\leq \frac{\sqrt{d(u,v)}}{1+\tau\sqrt{d(u,v)}}
$$
or
$$
\tau-\frac{1}{\sqrt{d(Tu,Tv)}}\leq-\frac{1}{\sqrt{d(u,v)}}
$$
and the contractivity condition of the Theorem 3 is satisfied with the function $\varphi(t)=-\frac{1}{\sqrt{t}}$ which belongs to the class $\mathfrak{F}$.

Consequently, by Theorem 3, the operator $T$ has a unique fixed point in $P$. This means that Problem (1) has a unique non-negative solution in $C[0,1]$. This finishes the proof.

An interesting question from a practical point of view is that the solution of Problem (1) is positive. A sufficient condition for that solution is positive, appears in the following result:

\textbf{Theorem 7.} Let assumptions of Theorem 6 be valid. If the function $t^\sigma f(t,y)$ is non-decreasing respect to the variable $y$, then the solution of Problem (1) given by Theorem 6 is positive.

\textbf{\emph{Proof:}} In contrary case, we find $t^*\in (0,1)$ such that $u(t^*)$=0. Since $u(t)$ is a fixed point of the operator $T$ (see Theorem 6) this means that
$$
u(t)=\int\limits_0^1 G(t,s) f(s,u(s))ds\,\, for\,\, 0<t<1.
$$
Particularly,
$$
0=u(t^*)=\int\limits_0^1 G(T^*,s)f(s,u(s))ds.
$$
Since that $G$ and $f$ are non-negative functions, we infer that
$$
G(t^*,s)f(s,u(s))=0\,\,\,a.e.\,(s)\eqno (2)
$$
On the other hand, as $\lim\limits_{t\rightarrow0^+}f(t,0)=\infty$ for given $M>0$ there exists $\delta>0$ such that for $s\in (0,\delta)$ $f(s,0)>M$. Since $t^\sigma f (t,y)$ is increasing and $u(t)\geq 0$,
$$
s^\sigma f(s,u(s))\geq s^\sigma f(s,0)\geq s^\sigma M\,\,\,for\,\,\, s\in (0,\delta)
$$
and, therefore, $f(s,u(s))\geq M$ for $s\in(0,\delta)$ and $f(s,u(s))\neq 0 \,\,a.e.\,\,(s)$. But this is a contradiction since $G(t^*,s)$ is a function of rational type in the variable $s$ and, consequently, $G(t^*,s)\neq 0 \,\,a.e.\,\,(s)$. Therefore, $u(t)>0$ for $t\in (0,1)$. This finishes the proof.

\section{Acknowledgement}
The present work is partially supported by the project MTM-2013-44357-P.

\end{document}